\documentclass[final]{elsart}

\usepackage{xspace,bm}
\usepackage{amssymb}
\usepackage{array}
\usepackage{epsfig,psfrag,graphics}
\usepackage{version}

\makeatletter

\theoremstyle{plain}

\newtheorem{theorem}{Result}

\renewcommand {\vec}[1]{{\bm{#1}}}
\renewcommand {\e}[1]{\mbox{e}^{#1}}
\def\Rset{{\mathbb{R}}}
\def\Nset{{\mathbb{N}}}
\def\A{{\mathcal A}}
\def\D{{\mathcal D}}
\def\S{{\mathcal S}}
\def\C{{\mathcal C}}
\newcommand {\ie}{{\it i.e.}\xspace}
\newcommand {\eg}{{\it e.g.}\xspace}
\begin{document}

\begin{frontmatter}

\title{Some extensions of the uncertainty principle}
\author[GIPSA]{Steeve     Zozor},     \author[IFLP]{Mariela     Portesi}     and
\author[UMLV]{Christophe Vignat}

\address[GIPSA]{GIPSA-Lab, D\'epartement Images et Signal\\
  961 Rue  de  la    Houille    Blanche, B.P. 46\\
  38420 Saint Martin d'H\`eres Cedex, France\\
  E-mail: {\tt steeve.zozor@gipsa-lab.inpg.fr}}

\address[IFLP]{Instituto de F\'isica La Plata (CONICET), and  \\
  Departamento de  F\'isica, Facultad de Ciencias  Exactas, Universidad Nacional
  de La  Plata\\
  C.C.\ 67, 1900 La  Plata, Argentina\\
  E-mail: {\tt portesi@fisica.unlp.edu.ar}}

\address[UMLV]{Laboratoire d'Informatique de l'Institut Gaspard Monge,\\
  Equipe  SYSCOM\\Universit\'e  de  Marne-la-Vall\'ee,  77454  Marne-la-Vall\'ee
  Cedex 2, France }

\begin{abstract}
  We study the formulation of  the uncertainty principle in quantum mechanics in
  terms  of  entropic  inequalities,   extending  results  recently  derived  by
  Bialynicki-Birula~\cite{Bia06} and Zozor  {\it et al.}~\cite{ZozVig07}.  Those
  inequalities  can   be  considered   as  generalizations  of   the  Heisenberg
  uncertainty principle,  since they  measure the mutual  uncertainty of  a wave
  function and its Fourier  transform through their associated R\'enyi entropies
  with conjugated indices.  We consider here the general case where the entropic
  indices are  not conjugated, in both  cases where the state  space is discrete
  and  continuous:  we  discuss  the  existence  of  an  uncertainty  inequality
  depending on the location of the  entropic indices $\alpha$ and $\beta$ in the
  plane $\left(\alpha, \beta \right)$.  Our  results explain and extend a recent
  study by  Luis~\cite{Lui07}, where states with quantum  fluctuations below the
  Gaussian case are discussed at the single point $(2,2)$.
\end{abstract}
\begin{keyword}
  Entropic   uncertainty  relation,   R\'enyi  entropy,  non-conjugated
  indices

\PACS 03.65.Ca, 03.65.Ta, 03.65.Db, 89.70.Cf
\end{keyword}
\end{frontmatter}


\section{Introduction}
\label{intro:sec}

The Uncertainty  Principle (UP)  is such a  fundamental concept that  it focuses
great attention not only in quantum  physics but even in other areas ({\it e.g.}
signal or image processing).  In  quantum mechanics terms, the UP etablishes the
existence of an  irreducible lower bound for the uncertainty in  the result of a
simultaneous measurement of non-commuting observables. An alternative expression
is  that the  precision  with  which incompatible  physical  observables can  be
prepared is limited by an upper bound.  Quantitatively, the UP can be given by a
relation of the form $U(A,B; \Psi)\geq \mathcal{B}(A,B)$, where $U$ measures the
uncertainty  in the  simultaneous  preparation  or measurement  of  the pair  of
operators  $A$  and $B$  when  the  quantum system  is  in  state $\Psi$,  while
$\mathcal{B}$  is a  state-independent bound.   The quantity  $U$ should  take a
fixed  minimum value  if  and only  if $\Psi$  is  a common  eigenstate of  both
operators.  It is interesting to note, as remarked by Deutsch~\cite{Deu83}, that
the quantitative  formulation of the UP through  generalizations of Heisenberg's
inequality  to  an  arbitrary  pair  of non-commuting  observables  (other  than
position and momentum)  may present some drawbacks, as  the expectation value of
the commutator between both operators (unless it is a $c$-number) does depend on
the  current  state  of  the  system.   Many authors  have  contributed  to  the
formulation of alternative  quantitative expressions for the UP,  among which is
the use of {\it entropic}  measures for the uncertainty, inspired by information
theory (see, for instance,  \cite{PorPla96} and references therein).  We address
here the  search for  lower bounds  of uncertainty relations  given in  terms of
generalized entropies of the R\'enyi form.

The paper is  organized as follows.  In section  \ref{recalls:sec}, we give some
brief recalls  on Fourier transform  in the context  of quantum physics,  and on
R\'enyi entropy, which can be viewed  as an extension of Shannon entropy.  Also,
we  summarize properties  of  entropic uncertainty  inequalities  for which  the
R\'enyi entropies  have conjugated  indices. Section \ref{arbitrary:sec}  is the
core  of this  work  :  its aim  is  to extend  the  usual entropic  uncertainty
inequalities  using  R\'enyi  entropies  with  arbitrary indices.   We  show  in
particular that  an uncertainty principle in  such a form does  not always exist
for  arbitrary pairs of  indices. In  section \ref{sec:conclusions},  we present
some conclusions,  while detailed proofs  of our main  results are given  in the
appendices at the end of the paper.


\section{Previous results on entropic uncertainty relations}
\label{recalls:sec}

We consider  a $d$-dimensional operator  $A$ and we  assume that the state  of a
system  is   described  by   the  wavefunction  $\Psi$.    Let  us   now  denote
$\widehat{\Psi}$ the Fourier transform of  $\Psi$, assumed to describe the state
of  operator  $\widetilde{A}$.  In  the  following, we  will  refer  to $A$  and
$\widetilde{A}$   as   conjugate  operators,   {\it   i.e.}  operators   having
wavefunctions that are linked by a Fourier transformation. As an example, we may
consider  the position  and momentum  of  a particle,  or the  position and  its
angular momentum.

The definition  of Fourier transform depends  on the state  space considered: if
the state space is discrete taking values  on the alphabet $\A^d = \{\ldots \: ,
\: 0 \:  , 1 \: , \ldots  \}^d$, finite or not, the Fourier  transform of $\Psi$
takes the form
\begin{equation}
\widehat{\Psi}(\vec{x}) = (2 \pi)^{- \frac{d}{2}} \sum_{\vec{a} \in \A^d}
\Psi(\vec{a}) \e{- \imath \vec{a}^t \vec{x}},
\label{rFT:eq}
\end{equation}
where $\vec{x} \in [0  \: ; \: 2 \pi)^d$.  In the  particular case of a discrete
and finite alphabet  of size $n$, $\A = \{0  \: , \: \ldots \:  , \: n-1\}$, one
can also consider a discrete Fourier transform
\begin{equation}
\widehat{\Psi}_{\mathrm{df}}(\vec{k}) = n^{- \frac{d}{2}} \sum_{\vec{a} \in
\A^d} \Psi(\vec{a}) \e{- \imath 2 \pi \vec{a}^t \vec{k}/n},
\label{DFT:eq}
\end{equation}
where $\vec{k} \in \A^d$.  When dealing with periodic quantities such as angular
momentum,   Eq.~(\ref{rFT:eq})  is   to  be   considered  \cite{Bia06,BiaMad85};
Eq.~(\ref{DFT:eq}) describes  a system where  the state and its  conjugate state
are both finite and discrete.

Finally, in the continuous case, the Fourier transform
takes the form
\begin{equation}
\widehat{\Psi}(\vec{x}) =  (2 \pi)^{- \frac{d}{2}}  \int_{\Rset^d} \Psi(\vec{u})
\e{- \imath \vec{u}^t \vec{x}} \, d\vec{u}.
\label{FT:eq}
\end{equation}

\

Let us consider now the R\'enyi $\lambda$-entropy associated to the operator $A$
when the physical system is in a state described by the wavefunction $\Psi$,
\begin{equation}
H_\lambda(A) = \frac{2 \lambda}{1-\lambda} \ln \|\Psi\|_{2 \lambda},
\label{Renyi:eq}
\end{equation}
for any  real positive  $\lambda \neq 1$,  where $\|.\|_p$ denotes  the standard
$p$-norm   of    a   function:   $\|\Psi\|_p    =   \left(   \int_{\mathbb{R}^d}
  |\Psi(\vec{u})|^p d\vec{u}  \right)^{1/p}.$ When $\lambda  \rightarrow 1$, the
R\'enyi  entropy $H_\lambda$  converges to the usual Shannon entropy
\begin{equation}
H_1(A) = - \int_{\mathbb{R}^d}  |\Psi(\vec{u})|^2 \ln |\Psi(\vec{u})|^2 d\vec{u}.
\end{equation}
When $\lambda  = 0$ one gets  $H_0(A) = \ln \mu(\{  x, \Psi(x) \ne  0 \})$ where
$\mu(.)$ is the  Lebesgue measure (\ie $H_0$ equals the  logarithm of the volume
of the  support of $\Psi$).  At the  opposite, when $\lambda \to  + \infty$, one
has $H_\infty = - \ln \|\Psi\|_{\infty} = - \ln \sup_{x \in \Rset^d} |\Psi(x)|$.

The  extension to a  discrete state  space is  straightforward by  replacing the
integral by  a discrete sum.  Equivalently, the  R\'enyi $\lambda$-entropy power
associated  to $A$ can  be defined  as\footnote{In \cite[p.   499]{CovTho91} the
  R\'enyi entropy  power is  defined as $  \exp(2 H_\lambda/d)$:  our definition
  does not affect the content of the paper and is in concordance with almost all
  cited papers. In  the Shannon case, the definition  $ \frac{1}{2 \pi \mbox{e}}
  \exp(2 H_\lambda/d)$ can also be found: with this last definition, the entropy
  power and the variance coincide in the Gaussian situation.}
\begin{equation}
N_\lambda(A) = \exp\left( \frac{1}{d} \, H_\lambda(A) \right).
\label{RenyiPower:eq}
\end{equation}
In the one-dimensional  finite discrete case, this entropy  power is the inverse
of the {\it certainty} measure  $M_{\lambda-1}$ defined by Maassen and Uffink in
Ref.~\cite{MaaUff88}.

\begin{figure}[ht]
\psfrag{a}{$\lambda$}
\psfrag{Na}{$N_\lambda$}
\begin{minipage}{.45\textwidth}
\centerline{\includegraphics[height=4cm]{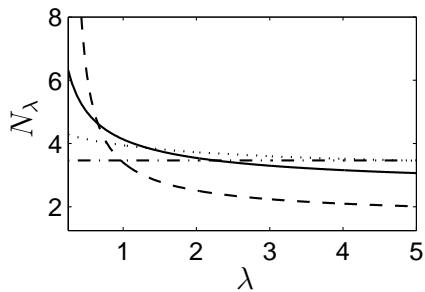}}
\end{minipage}
\begin{minipage}{.55\textwidth}
  \caption{Behavior of  $N_\lambda $ versus $\lambda$  for various wavefunctions
    (identity   covariance  matrix,   $d=1$)~:  Gaussian   case   (solid  line),
    Student-$t$ case  $\Psi(\vec{x}) \propto (  1 + \vec{x}^t \vec{x}  / (\nu-2)
    )^{-\frac{d+\nu}{4}}$   for   $\nu=3$   (dashed  line),   Student-$r$   case
    $\Psi(\vec{x}) \propto ( 1 - \vec{x}^t \vec{x} / (\nu+2))^{\frac{\nu-d}{4}}$
    for $\nu=3$  (dotted line) and Uniform  on the sphere,  \ie Student-$r$ with
    $\nu=d$ (dash-dotted line).}
\label{Nlambda:fig}
\end{minipage}
\end{figure}

Among the  properties of the one-parameter  family of R\'enyi  entropies, let us
mention that for arbitrary fixed  $A$ (or $\Psi$), $H_\lambda$ is non increasing
versus $\lambda$, and hence $N_\lambda$ is also non increasing against $\lambda$
\cite[th.  192]{HarLit52}.  This  property can easily be shown  by computing the
derivative against $\lambda$, the derivative of the factor of $1/(1-\lambda)^2$,
and  using the  Cauchy-Schwartz inequality.   How $N_\lambda$  decreases against
$\lambda$  is then  intimately linked  to $\Psi$:  in the  particular case  of a
constant  $\Psi$  (on  a  bounded  support,  in  order  to  ensure  wavefunction
normalization),  $N_\lambda$   is  constant.   This  is   illustrated  in  figure
\ref{Nlambda:fig}  where  the entropy  power  $N_\lambda(A)$  is plotted  versus
$\lambda$ for  various observables $A$ with  corresponding wavefunctions $\Psi$.
This figure  illustrates also  the fact  that the maximal  entropy power  is not
given by the same $\Psi$ for any $\lambda$~: under covariance matrix constraint,
$N_\lambda$ is maximal in the  Gaussian case for $\lambda=1$, in the Student-$t$
case with $\nu$ degrees of freedom for $\lambda = 1 - \frac{2}{d+\nu}$ or in the
Student-$r$   case  with  $\nu$   degree  of   freedom  for   $\lambda  =   1  +
\frac{2}{\nu-d}$ \cite{CosHer03,VigHer04}.

\

Uncertainty  relations involving $A$  and $\tilde  A$ that  make use  of certain
combinations of  R\'enyi generalized entropies have already  been established in
Ref.~\cite{Bia06,ZozVig07}  in the  continuous--continuous, discrete--continuous
(periodic)   and   discrete--discrete   cases.    In  terms   of   the   R\'enyi
$\lambda$-entropy power, these uncertainty relations read
\begin{equation}
N_{\frac{p}{2}}(A) \, N_{\frac{q}{2}}(\widetilde{A}) \ge C_{p,q},
\label{UP:eq}
\end{equation}
where the  entropic parameters $\frac{p}{2}$ and $\frac{q}{2}$  have been chosen
{\it  conjugated}\footnote{Rigorously speaking,  $p$ and  $q$ are  conjugated if
  $1/p+1/q=1$. Thus  here, $p/2$ and $q/2$  are not conjugated;  however, in the
  rest of  the paper, and  for the  sake of simplicity,  we will call  $p/2$ and
  $q/2$     such    that     $1/p+1/q=1$    ``conjugated''     indices.},    \ie
$\frac{1}{p}+\frac{1}{q} =  1$, and  with $p \ge  1$ (infinite $q$  when $p=1$).
The constant lower bound is expressed as
\begin{equation}
C_{p,q} = \left\{\begin{array}{ll}
n & \mbox{ in the discrete--discrete case}\\
2 \pi p^{\frac{1}{p-2}} q^{\frac{1}{q-2}} & \mbox{ in the continuous--continuous case}\\
2 \pi & \mbox{ in the discrete--continuous case}.
\end{array}\right.
\label{Bound:eq}
\end{equation}
In  the  continuous-continuous  case,  we  extend  by  continuity  $C_{p,q}$  as
$C_{1,+\infty} =  2 \pi$  and $C_{2,2}  = \mbox{e} \pi$.   In the  discrete case
(finite  or not), this  kind of  uncertainty relation  is usually  exhibited for
one-dimensional  quantities, but  it  extends trivially  to the  $d$-dimensional
situation by  one-to-one mapping between $\Nset^d$  and $\Nset$; it  is based on
the  Hausdorff inequality \cite[th.   IX.8, p.11]{ReeSim75},  which is  itself a
consequence  of   the  Riesz--Thorin  interpolation   theorem  \cite[th.   IX.17
p.27]{ReeSim75}.  These theorems also give a  lower bound for the product of the
entropy   powers  $N_{\frac{p}{2}}(A)  N_{\frac{q}{2}}(\widetilde{A})$   in  the
continuous--continuous case, but this bound is not sharp.  However, this case is
addressed  via the  Beckner  relation given  in  Ref.~\cite{Bec75:02}, as  shown
in~\cite{Bia06,ZozVig07,BiaMyc75}.    Finally,  we   remark  that   the  product
$N_{\frac{p}{2}}(A) \, N_{\frac{q}{2}}(\widetilde{A})$ is scale invariant in the
sense that it  is independent of invertible matrix  $M$ when replacing $\Psi(x)$
by $|M|^{-\frac{d}{2}} \Psi(M^{-1} x)$.

It is  worth stressing that inequality~(\ref{UP:eq}) provides  a universal (\ie,
independent  of the state  of the  system) lower  bound for  the product  of two
measures of uncertainty  that quantify the missing information  related with the
measurement or preparation of the  system with operators $A$ and $\widetilde A$.
We note that  in the continuous--continuous context, equality  is reached if and
only if  $\Psi$ is  a Gaussian wave  function; in  the case of  discrete states,
equality is attained if and only if $\Psi$ coincides with a Kronecker indicator,
or  with a  constant  (by conjugation,  and  provided that  the  space state  is
finite).  Moreover, maximization  of $(H_1(A)+H_1(\widetilde{A}))/d$ ($p/2 = q/2
=  1$  conjugated) without  constraint  has  been  suggested as  an  interesting
counterpart to the  usual approach of entropy maximization  under constraint for
the derivation of the wavefunction associated to atomic systems \cite{GadBen85}.

It is important  to note that all uncertainty  relations (\ref{UP:eq}) deal with
R\'enyi  entropies  with  {\it   conjugated}  indices.   Indeed,  the  Hausdorff
inequality,  as   a  consequence  of   the  Parseval  relation   expressing  the
conservation  by Fourier  transformation of  the  Euclidean norm,  and from  the
trivial  relation  $\|\Psi\|_\infty \le  \|\Psi\|_1$,  involves only  conjugated
indices, as recalled \eg in  Ref.~\cite{MaaUff88}.  In the following section, we
study  extensions of  these uncertainty  inequalities to  the case  of arbitrary
pairs of {\it non-conjugated} indices.   Such a situation has been considered by
Luis in  \cite{Lui07} for  a particular combination  of indices to  describe the
uncertainty product of exponential states.   Moreover, one may hope to gain more
flexibility  by using non-conjugated  information measures  in both  spatial and
momentum  domains.  As  a supplementary  motivation,  one may  wish to  quantify
uncertainty  for   conjugate  observables   using  the  same   entropic  measure
$N_\lambda$:  inequality (\ref{UP:eq})  then  holds only  when $p=q=2$  (Shannon
case).


\section{Uncertainty relations with entropies of arbitrary indices}
\label{arbitrary:sec}

We  present  here an  extension  of  entropic  formulations of  the  Uncertainty
Principle for conjugate observables, to the general situation of arbitrary pairs
of entropic indices. For this purpose, we discuss separately the three different
cases corresponding to operators having discrete or continuous spectrum.


\subsection{Discrete--discrete case}

Inequality  (\ref{UP:eq})  appears in  Ref.~\cite{MaaUff88}  in  the context  of
finite (discrete)  states and $d=1$,  but an extension to  arbitrary nonnegative
indices $\alpha$  and $\beta$ can also be  provided \cite{MaaUff88,Lar90}, which
reads with the notation adopted here
\begin{equation}
N_\alpha(A) N_\beta(\widetilde{A}) \ge \left(\frac{2 n}{n+1}\right)^2
\label{UPMaassen:eq}
\end{equation}
This  uncertainty relation  is first  proved in  the case  $\alpha$  and $\beta$
infinite and then, for any pair of indices $\alpha$ and $\beta$ as a consequence
of the  positivity and the  non-increasing property of $N_\lambda$  in $\lambda$
(see  above and  \cite[th.   16]{HarLit52}).  This  inequality was  rediscovered
recently   by  Luis   in   the  particular   case   of  non-conjugated   indices
$\alpha=\beta=2$  \cite[eq.    2.8]{Lui07}.   Again,  equality   is  reached  in
(\ref{UPMaassen:eq})  when  $\Psi$  coincides  with a  Kronecker  indicator  and
$\widehat{\Psi}$ is constant, or ``conjugately'', for a constant $\Psi$.


\subsection{Continuous--continuous case}


\subsubsection{A result by Luis}


In Ref.~\cite{Lui07}, Luis claims that, to his best knowledge, there is no known
continuous  counterpart  of  the   finite  discrete  case.   Furthermore,  after
mentioning  that  only  Gaussian  wavefunctions  saturate  both  the  Heisenberg
uncertainty relation  and the conjugated-indices  entropic uncertainty relation
(\ref{UP:eq}), Luis  considers the product  of entropy powers in  the particular
case  $\alpha  =  \beta  =  2$.  Indeed, he  compares  this  product,  for  some
particularly  distributed states,  to  the value  obtained  for Gaussian  states
$\Psi_G$  associated to  operator  $G$.  He  remarks  that this  product is  not
saturated by  Gaussians since,  for one-dimensional exponential  states $\Psi_E$
(resp. operator $E$), he finds that
\begin{equation}
N_2(E)  N_2(\widetilde{E}) =
\frac{8 \pi}{5} < 2 \pi = N_2(G) N_2(\widetilde{G}).
\end{equation}
In what  follows, we  explain this  result and deduce  extended versions  of the
entropic uncertainty relations for arbitrary indices.


\subsubsection{Uncertainty relations for arbitrary indices and domain of existence}

\begin{figure}[htbp]
\centerline{\includegraphics{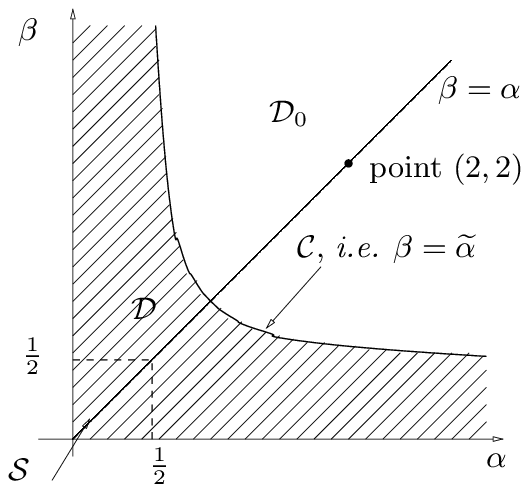}}
\caption{Sets $\D_0$ (blank  region), $\D$ (shaded region), $\S$  (square in the
  shaded region) and  $\C$ (solid line) in the  $(\alpha,\beta)-$plane, as given
  by Eq.~(\ref{domains:eq}). Notice that, by definition, $\C \subset \D$}
\label{domains:fig}
\end{figure}

For any real $\alpha \ge \frac{1}{2}$, let us introduce the notation
\begin{equation}
\tilde{\alpha} = \frac{\alpha}{2 \alpha-1}
\label{tilde:eq}
\end{equation}
so that  $\alpha$ and $\tilde{\alpha}$  are conjugated indices, \ie  $1/\alpha +
1/\tilde{\alpha}=2$ (for $\alpha =  1/2$, $\tilde{\alpha}$ is infinite), and the
function $B: [1/2 \: ; \: +\infty) \mapsto [2\pi \: ; \: e\pi]$
\begin{equation}
B(\alpha) = \pi \alpha^{\frac{1}{2(\alpha-1)}} \tilde{\alpha}^{\frac{1}{2(\tilde{\alpha}-1)}}
\label{B:eq}
\end{equation}
When $\alpha \to \frac{1}{2}^+$, $B(\alpha) \to 2 \pi$; hence, by continuity, we
set $B(1/2) = 2 \pi$. By continuity, we also set $B(1) = \mbox{e} \pi$.  We also
define the following sets on the plane
\begin{equation}\left\{\begin{array}{lll}
\D_0 & = & \left\{ (\alpha,\beta) \: : \: \alpha > \frac{1}{2} \mbox{ and } \beta
        > \tilde{\alpha} \right\}\\
\S & = & \left[0 \: ; \: \frac{1}{2} \right)^2\\
\D & = & \Rset_+^2 \backslash \D_0\\
\C & = & \left\{(\alpha,\beta) : \alpha \ge \frac{1}{2}, \beta=\tilde{\alpha} \right\}.
\end{array}\right.
\label{domains:eq}
\end{equation}
which are represented in Fig.~\ref{domains:fig}. The solid line represents curve
$\C$  for  which $\alpha$  and  $\beta$ are  conjugated  and  the shaded  region
represents $\D$.

We recall that most of the previous results on uncertainty relations in terms on
$\lambda$-R\'enyi  entropies  refer  to  inequalities for  pairs  of  conjugated
entropic    indices    \ie    located    on   curve    $\cal{C}$    (see    Eqs.
(\ref{UP:eq})-(\ref{Bound:eq}) above).  We first  remark that, as for conjugated
indices, the product $N_\alpha(A) N_\beta(\tilde{A})$ is scale invariant for any
pair  of indices  $(\alpha,\beta) \in  \Rset_+^2$. We  now introduce  some novel
results  for  {\em arbitrary  pairs  of indices},  depending  on  the region  of
$\Rset_+^2$ where they lie.

\begin{theorem}
  For  any pair $(\alpha,\beta)  \in \D$  and for  conjugate operators  $A$ and
  $\widetilde{A}$, there exists an uncertainty principle under the form
\begin{equation}
N_\alpha(A) N_\beta(\widetilde{A}) \ge B_{\alpha,\beta},
\label{generalUP_C:eq}
\end{equation}
where
\begin{equation}
B_{\alpha,\beta} = \left\{\begin{array}{lll}
B(\alpha) & \mbox{ in } & (\D \backslash \S) \cap \{ (\alpha,\beta) : \alpha \ge \beta\}\\
B(\beta) & \mbox{ in } & (\D \backslash \S) \cap \{ (\alpha,\beta) : \beta \ge \alpha\}\\
2 \pi & \mbox{ in } & \S.
\end{array}\right.
\label{generalBound_C:eq}
\end{equation}
\label{thCCUP:th}
\end{theorem}

The proof of this result is given in appendix \ref{result1:app}. Except on $\C$,
the bound $B_{\alpha,\beta}$ is probably not sharp and we have not determined if
such an uncertainty  saturates for Gaussians or not.   A direct calculation with
Gaussians    shows    that    $N_\alpha(G)    N_\beta(\widetilde{G})    =    \pi
\alpha^{\frac{1}{2 (\alpha-1)}} \beta^{\frac{1}{2 (\beta-1)}}$ which is strictly
higher  than  $B_{\alpha,\beta}$  in  $\D  \backslash \C$:  thus,  either  bound
(\ref{generalBound_C:eq}) is  not sharp in  $\D \backslash \C$, or  Gaussians do
not saturate (\ref{generalUP_C:eq})--(\ref{generalBound_C:eq}) in $\D \backslash
\C$, or  both. This  point remains to  be solved.   Note however that  the point
(0,0) is  degenerate since $N_0(A)$ measures  the volume of the  support of $A$.
It is  well known  that if  $\Psi$ is defined  on a  finite volume  support, the
support  of   its  Fourier  transform  has  infinite   measure.   Hence  $N_0(A)
N_0(\widetilde{A})     =     +      \infty$     which     trivially     fulfills
(\ref{generalUP_C:eq})-(\ref{generalBound_C:eq}).   By continuity  of $N_\alpha$
in $\alpha$,  this remark suggests  that the bound  (\ref{generalBound_C:eq}) is
not sharp, at least in $\S$.

Our  second   result  concerns  the  non-existence  of   a  generalized  entropy
formulation for the uncertainty principle in the area $\D_0$.
\begin{theorem}
  For any  pair $(\alpha,\beta) \in \D_0$  and for conjugate  operators $A$ and
  $\widetilde{A}$,  the  positive product  of  the  entropy powers  $N_\alpha(A)
  N_\beta(\widetilde{A})$ can be arbitrarily small.  In other words, no entropic
  uncertainty principle exists in $\D_0$.
\label{D0:th}
\end{theorem}

The   proof  and   illustrations  of   this   result  are   given  in   appendix
\ref{result2:app}.   Figure~\ref{UPDomainsBounds:fig} schematizes  the preceding
results: entropic uncertainty relations exist in $\D$ but not in $\D_0$.

\begin{figure}[htbp]
\centerline{\includegraphics{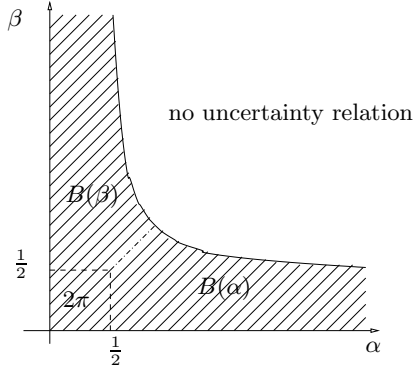}}
\caption{Entropic  uncertainty relations  exist for  $(\alpha,\beta) \in  \D$ as
  depicted by  the shaded  area, with the  corresponding bound.   Conversely, in
  $\D_0$  the  positive  product  $N_\alpha(A)  N_\beta(\widetilde{A})$  can  be
  arbitrarily small.}
\label{UPDomainsBounds:fig}
\end{figure}

In the case studied by Luis \cite{Lui07}, $(\alpha,\beta) = (2,2) \in \D_0$ (see
Fig. \ref{domains:fig}) so  that one can find wave  functions (or operators) for
which $N_2(A)  N_2(\widetilde{A})$ can be  arbitrarily small.  This is  thus not
surprising that Luis finds a wavefunction  for which this entropy power is lower
than in  the Gaussian  case: he  considers in fact  the special  one dimensional
Student-t  case  with  $\nu=3$  of Eq.~(\ref{Student-t_hat:eq})  below;  varying
parameter $\nu$, one  can describe a richer family  of distributions that allows
to  prove result~\ref{D0:th}.   Other cases  having equal  entropic  indices are
worth studying in more details.   These cases correspond to points located along
the  line  that  bisects  $\Rset_+^2$  (see  Fig.   \ref{domains:fig})  and,  as
mentioned, the situation is very different whether the point lies inside $\D$ or
inside $\D_0$.

\begin{figure}[ht]
\psfrag{a}{$\alpha$}
\psfrag{NaNa}{\hspace{-7.5mm} \footnotesize $N_\alpha(A) N_\alpha(\widetilde{A})$}
\begin{minipage}{.45\textwidth}
\centerline{\includegraphics[height=4cm]{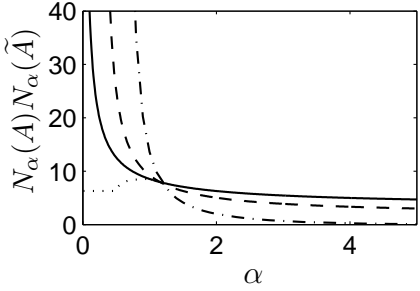}}
\end{minipage}
\begin{minipage}{.55\textwidth}
  \caption{Behavior of $N_\alpha(A) N_\alpha(\widetilde{A})$ versus $\alpha$ for
    various wavefunctions ($d=1$)~: Gaussian case (solid line), Student-$t$ case
    $\Psi(\vec{x})     \propto     (     1     +     \vec{x}^t     \vec{x}     /
    (\nu-2))^{-\frac{d+\nu}{4}}$ for $\nu=3$  (dashed line) and Student-$t$ case
    for $\nu=.8$ (dash-dotted  line).  The dotted line depicts  the lower bounds
    (\ref{generalBound_C:eq}) that only exist when $\alpha \le 1$.}
\label{NaNa_a:fig}
\end{minipage}
\end{figure}

Figure  \ref{NaNa_a:fig}  depicts  the  behavior  of  the  product  $N_\alpha(A)
N_\alpha(\widetilde{A})$  versus  $\alpha$, for  $\alpha  >0$,  in the  Gaussian
context,  in  the  Student-$t$/Laplace  context  of  Luis,  and  in  a  specific
Student-$t$ case used  in the proof of result  \ref{D0:th}. This illustration is
motivated  by  the interest  in  using the  same  entropic  measure to  describe
uncertainty for conjugate observables.  It clearly seen that, when $\alpha > 1$,
no uncertainty principle holds.   Moreover, for many operators $A$, $N_\alpha(A)
N_\alpha(\widetilde{A})$  can  be  below  the  product  for  the  Gaussian  case
$N_\alpha(G) N_\alpha(\widetilde{G})$, for a wide range of index $\alpha$.

Concerning our second result for any pair of R\'enyi indices $(\alpha,\beta) \in
\D_0$  and   conjugate  operators  $A$  and   $\widetilde{A}$,  two  alternative
interpretations  arise: either  uncertainty principle  does not  apply  in those
cases,  which  contradicts  physical  intuition, or  the  R\'enyi  entropy-power
product is not suitable to quantify uncertainty for conjugate observables in the
area  $\D_0$.   The R\'enyi  as  well  as Shannon  entropy  related  to a  given
observable quantifies the amount of  missing information in the knowledge of the
properties of a  system in connection with a measurement  of that observable. It
is  in this  sense  that  one can  consider  $H_\lambda$ as  a  measure of  {\it
  uncertainty}. In order  not to contradict the common  principle of uncertainty
in quantum  physics for pairs of  non-commuting observables, one  can argue that
for $\alpha$ fixed and a  given quantum observable, the $\beta$-R\'enyi entropy,
when  $\beta$ is  such that  $(\alpha,\beta) \in  \D_0$, is  not well  suited to
describe  the  lack  of  knowledge   about  the  conjugate  observable.   It  is
interesting to note that for  $\alpha=\beta=2$, the entropy power product is the
Onicescu  measure  used  to  quantify  complexity or  disequilibrium  in  atomic
systems\cite{ChaMon05,SenPan07}.   However, it  has been  suggested to  use this
product  $N_2(A)  N_2(\widetilde{A})$, divided  by  the  product  of the  Fisher
informations  related to  $A$  and  $\widetilde{A}$. This  can  suggest that  in
itself, $N_2(A) N_2(\widetilde{A})$  is not a reasonable measure  to account for
uncertainty for conjugate observables.


\subsection{Discrete--continuous case}

Our last result concerns the  discrete--continuous case, and can be expressed as
follows.
\begin{theorem}
  For   any  pair   $(\alpha,\beta)  \in   \D$   and  for   conjugate  $A$   and
  $\widetilde{A}$, there exists an uncertainty principle under the form
\begin{equation}
N_\alpha(A) N_\beta(\widetilde{A}) \ge 2 \pi.
\end{equation}
\label{thCDUP:th}
\end{theorem}

The proof  is similar to that  of result~\ref{thCCUP:th}.  The  bound comes from
(\ref{Bound:eq}).   We  remark  that  here,  the  lower  bound  of  $N_\alpha(A)
N_\beta(\widetilde{A})$ is sharp and  is attained when $\Psi(\vec{k})$ coincides
with a Kronecker indicator.

What happens in the $\D_0$ domain remains to be solved.


\section{Conclusions}
\label{sec:conclusions}
We have addressed some fundamental questions related with the formulation of the
Uncertainty Principle for  pairs of conjugate operators, like  \eg position and
momentum,  in  entropic  terms.   Our  study  extends  the  set  of  uncertainty
inequalities as derived by Bialynicki-Birula  \cite{Bia06} to the case where the
entropic  indices are  not  conjugated.   Our main  findings  are summarized  in
results~\ref{thCCUP:th},  \ref{D0:th}  and  \ref{thCDUP:th} that  establish  the
conditions  under which  an entropic  formulation of  the  Uncertainty Principle
makes sense and, if so, its lower bound.  We have addressed the cases where both
state space  and Fourier transformed  state space are respectively  (i) discrete
and discrete, (ii) continuous and continuous, and (iii) discrete and continuous.
The cases of equality in the  uncertainty relation considered (\ie, the state of
the system corresponding to  minimum uncertainty in the simultaneous measurement
of both  observables) are still undetermined  and will be the  object of further
research.

Summing up, our study establishes very general conditions for the formulation of
the  Uncertainty  Principle  in  entropic   terms  (as  an  alternative  to  the
Robertson--Schr\"odinger  formulation  in terms  of  variances),  making use  of
generalized  entropies  as  measures  of  uncertainty  for  the  preparation  or
measurement  of pairs  of quantum  observables in  a given  state of  a physical
system.  We believe that our analysis  sheds some light on previous related work
in the  field, and also  that it has  implications in the discussion  of quantum
behavior of physical systems, like quantum limits to precision measurements, for
instance.


\appendix

\section{Proof of Result 1}
\label{result1:app}

One  deals here  with entropic  uncertainty  products of  the form  $N_\alpha(A)
N_\beta(\tilde A)$  for arbitrary pairs  of indices $(\alpha,\beta) \in  \D$ and
for conjugate operators  $A$ and $\tilde A$ having  continuous spectra. We first
prove existence and  then evaluate lower bounds in  different regions inside set
$\cal D$.

\begin{figure}[htbp]
  \centerline{\includegraphics{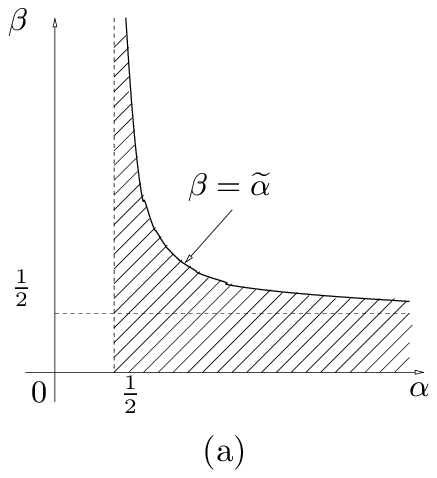}
    \hspace{15mm} \includegraphics{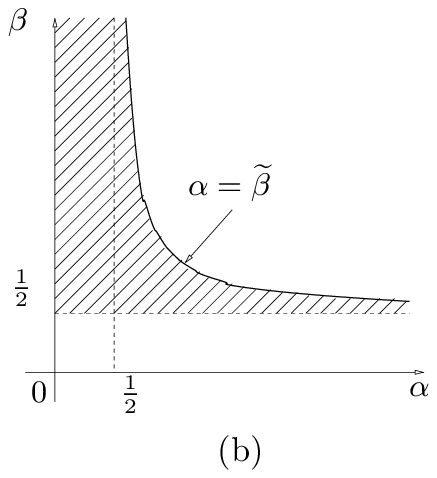}}
  \caption{(a) The dashed area represents  the area where for any fixed $\alpha$
    and for $\tilde{\alpha} \ge \beta \ge 0$, uncertainty relation (\ref{UP:eq})
    extends  to.  (b):  The  problem  is symmetric  by  exchanging $\alpha$  and
    $\beta$; the dashed  area represents then the ``conjugated''  area dashed in
    (a).}
\label{case_alpha_beta:fig}
\end{figure}

\begin{enumerate}
\item   From    (\ref{UP:eq})--(\ref{Bound:eq}),   setting   $p=2\alpha$   (then
  $q=2\tilde\alpha$), one has $N_\alpha(A) N_{\tilde{\alpha}}(\widetilde{A}) \ge
  B(\alpha)$. This proves the result for indices on curve $\cal C$.

\item In order  to prove the result  for indices in $\D \backslash  \C$ we first
  restrict the proof to $\D \backslash \S$ in two step:

  \begin{enumerate}
  \item  Fix  $\alpha \ge  \frac{1}{2}$:  since  the  R\'enyi entropy  power  is
    decreasing  in  $\beta$  \cite[th.    192]{HarLit52},  for  all  $\beta  \le
    \tilde{\alpha}$,   inequality    $N_\alpha(A)   N_\beta(\widetilde{A})   \ge
    B(\alpha)$     still    holds.     This     case    is     schematized    in
    Fig.~\ref{case_alpha_beta:fig}(a).\label{case_alpha:it}
  \item  The same  arguments apply  by  symmetry between  $\alpha$ and  $\beta$:
    $N_{\alpha}(A) N_{\beta}(\widetilde{A}) \ge B(\beta) $ provided that $\alpha \le
    \tilde{\beta}$.      This     case      is     schematized     in
    Fig.~\ref{case_alpha_beta:fig}(b).\label{case_beta:it}
  \end{enumerate}
  As a  conclusion, uncertainty relation  (\ref{generalUP_C:eq}) exists provided
  $\alpha \ge 1/2$ and $\beta \le \tilde{\alpha}$ (first case), or provided $\beta
  \ge  1/2$ and  $\alpha \le  \tilde{\beta}$ \ie  provided $(\alpha,\beta)  \in \D
  \backslash \S$.
\item   Let    us   consider   now   the   set    $\S$:   uncertainty   relation
  (\ref{generalUP_C:eq}) remains true  in the limit $\alpha \to  1/2$, and hence
  in  point $(1/2  ;  1/2)$.  Thus,  with  the same  argument of  non-increasing
  property  of  $N_\alpha$,  an  uncertainty  relation  again  exists  in  $\S$,
  $N_\alpha(A)       N_{\beta}(\widetilde{A})       \ge       N_{\frac{1}{2}}(A)
  N_{\frac{1}{2}}(\widetilde{A})              \ge             N_{\frac{1}{2}}(A)
  N_{\infty}(\widetilde{A})$.\label{case_S:it}
\end{enumerate}
As a conclusion, an uncertainty principle exists in the whole area $\D$.

In order to evaluate the  bound $B_{\alpha,\beta}$, we restrict our attention to
the   case   $\alpha   \ge    \beta$   since   trivially   $B_{\alpha,\beta}   =
B_{\beta,\alpha}$.
\begin{itemize}
\item   For  $\alpha   \ge  \frac{1}{2}$   and  $\beta   <   \frac{1}{2}$,  from
  (\ref{Bound:eq})    or    figure    \ref{case_alpha_beta:fig},    only    case
  (\ref{case_alpha:it})  is   to  be  considered.   Hence,  $B_{\alpha,\beta}  =
  B(\alpha)$.
\item For $\alpha \in \left(\frac{1}{2} \:  ; \: 1 \right]$, a study of function
  $B(\alpha)$ shows that  it increases in $\left( \frac{1}{2}  ; 1 \right]$ from
  value   $2  \pi$   to  $\pi   \mbox{e}$.    In  the   situation  $\alpha   \in
  \left(\frac{1}{2} \: ;  \: 1 \right]$ and $\beta  \ge \frac{1}{2}$, both cases
  (\ref{case_alpha:it})    and    (\ref{case_beta:it})    occur   (see    figure
  \ref{case_alpha_beta:fig})      but     since     $\beta      \le     \alpha$,
  $\max(B(\alpha),B(\beta))  =  B(\alpha)$.   Thus,  again  $B_{\alpha,\beta}  =
  B(\alpha)$.
\item When $\alpha > 1$, the study of  $B$ shows also that it decreases in $(1 ;
  + \infty)$ (from $\pi \mbox{e}$ to $2\pi$) and furthermore that for $\beta \le
  \tilde{\alpha}$,      $B(\beta)      \le      B(\alpha)$.      Here      again
  $\max(B(\alpha),B(\beta))   =   B(\alpha)$   and  thus   $B_{\alpha,\beta}   =
  B(\alpha)$.
\item   Finally,  when   $\alpha   <  1/2$,   from   case  \ref{case_S:it}   and
  $B\left(\frac{1}{2} \right)  = 2 \pi$, one  has $B_{\alpha,\beta} =  2 \pi$ in
  $\S$.
\end{itemize}


\section{Proof of Result 2}
\label{result2:app}

To  prove  this  result, it  is  sufficient  to  exhibit  an example  for  which
$N_\alpha(A)   N_\beta(\widetilde{A})$    can   be   arbitrarily    small   when
$(\alpha,\beta)  \in  \D_0$.   To  this  aim,  let  us  consider  the  following
Student-$t$ wavefunction
\begin{equation}
    \Psi(\vec{x}) =
\sqrt{\frac{\Gamma\left(\frac{d+\nu}{2}\right)}{\pi^{\frac{d}{2}}
\Gamma\left(\frac{\nu}{2}\right)}} \: \left( 1 + \vec{x}^t \vec{x}\right)^{- \,
\frac{d+\nu}{4}},
\label{Student-t:eq}
\end{equation}
where $\Gamma$ is the Gamma function and  $\nu > 0$ is a parameter called degree
of  freedom.   Its  Fourier   transform,  from  Refs.~\cite[eq.   5]{Lor54}  and
\cite[6.565-4]{GraRyz80}, reads
\begin{equation}
\widehat{\Psi}(\vec{x}) = \sqrt{\frac{2^{\frac{4-d-\nu}{2}}
\Gamma\left(\frac{d+\nu}{2}\right)}{\pi^{\frac{d}{2}}
\Gamma\left(\frac{\nu}{2}\right) \Gamma^2\left(\frac{d+\nu}{4}\right)}} \:
(\vec{x}^t \vec{x})^{\frac{\nu-d}{8}} \, K_{\frac{d-\nu}{4}}((\vec{x}^t
\vec{x})^{\frac{1}{2}}),
\label{Student-t_hat:eq}
\end{equation}
where $K_\mu$ is the modified Bessel function of the second kind of order $\mu$.
From Ref.~\cite[4.642 and 8.380-3]{GraRyz80}, the R\'enyi $\alpha$-entropy power
associated to $\Psi$ is expressed as
\begin{equation}
N_\alpha(A) = \sqrt{\pi} \, \left(\frac{\Gamma^\alpha
\left(\frac{d+\nu}{2}\right)}{\Gamma \left(\frac{\alpha
(d+\nu)}{2}\right)}\right)^{\frac{1}{d(1-\alpha)}} \: \left(\frac{\Gamma
\left(\frac{\alpha (d+\nu)-d}{2}\right)}{\Gamma^\alpha
\left(\frac{\nu}{2}\right)}\right)^{\frac{1}{d(1-\alpha)}} \mbox{ if } \alpha
\ne 1.
\label{NSt:eq}
\end{equation}
The case  $\alpha=1$ is obtained by continuity.   R\'enyi $\alpha$-entropy power
is   then  defined   provided  that   $   \nu  >   \max  \left(   0  ,   \frac{d
    (1-\alpha)}{\alpha}  \right)$, \ie  $\nu >  0$ if  $\alpha >  1$ and  $\nu >
\frac{d (1-\alpha)}{\alpha}$ otherwise.

Similarly, the R\'enyi $\beta$-entropy power associated to $\widehat{\Psi}$ is
expressed as
\begin{eqnarray}
N_\beta(\widetilde{A}) & = & \sqrt{\pi} \, \left( \frac{2^{\frac{(4-d-\nu) \beta
+ 2}{2 d}} \, \Gamma \left(\frac{d+\nu}{2}\right)}{\Gamma^2
\left(\frac{d+\nu}{4}\right)} \right)^{\frac{\beta}{d (1-\beta)}} \nonumber\\
& & \times \left( \frac{\displaystyle \int_0^{+\infty} 
r^{d-1+\frac{\beta(\nu-d)}{2}} K_{\frac{d-\nu}{4}}^{2\beta}(r) \,
dr}{\Gamma^\beta \left(\frac{\nu}{2}\right)} \right)^{\frac{1}{d(1-\beta)}}
\mbox{ if } \beta \ne 1,
\label{NSttilde:eq}
\end{eqnarray}
the case  $\beta=1$ being  obtained by continuity.   From the properties  of the
Gamma    function,   from   \cite[9.6.8    and   9.6.9]{AbrSte70}    and   since
$K_\mu=K_{-\mu}$, the R\'enyi $\beta$-entropy  power exists provided that $\nu >
\max \left( 0 , \frac{d (\beta-1)}{\beta} \right)$, \ie $\nu > 0$ if $\beta < 1$
and $\nu > \frac{d (\beta-1)}{\beta}$ otherwise.

\begin{figure}[htbp]
\psfrag{nu}{$\nu$}
\psfrag{NaNbP}{\hspace{-9mm}\footnotesize $N_\alpha, N_\beta$ and $N_\alpha N_\beta$}
  \centerline{\includegraphics[width=6.5cm]{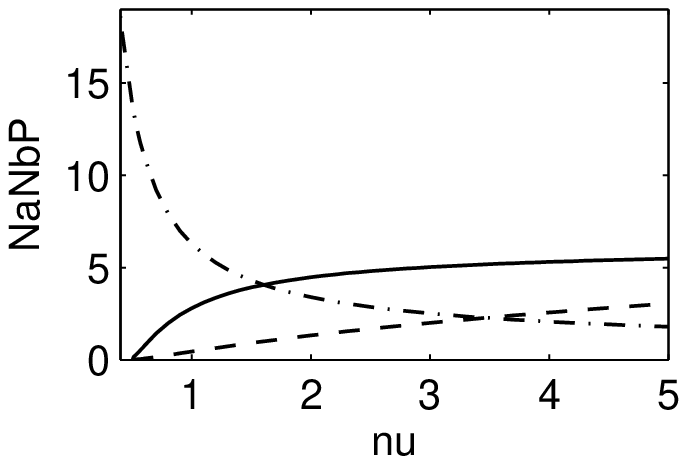}
    \hspace{5mm} \includegraphics[width=6.5cm]{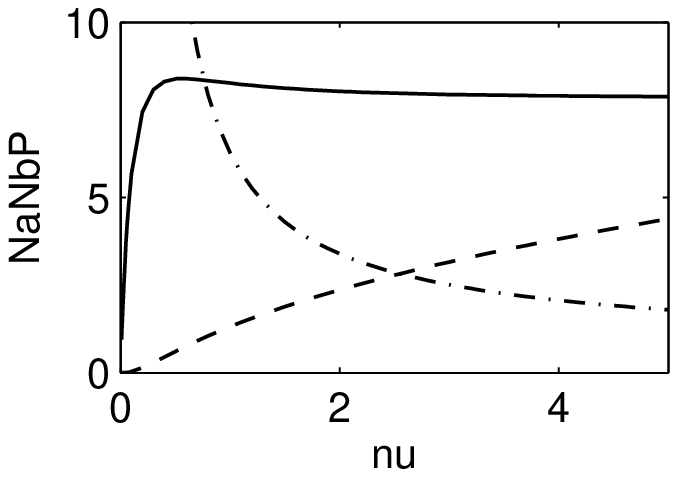}}
  \caption{Behavior       of       $N_\alpha(A)$       (dash-dotted       line),
    $N_\beta(\widetilde{A})$  (dashed  line),  and  of the  product  $N_\alpha(A)
    N_\beta(\widetilde{A})$ (solid  line) versus  $\nu$ in the  Student-$t$ case
    described  in  the  proof  of   result  \ref{D0:th}  (see  text).  In  these
    illustrations $d=1$,  $(\alpha,\beta) = (2,2)$ (left)  and $(\alpha,\beta) =
    (2,3/4)$ (right)  are both in $\D_0$. It  confirm that, when $\beta  > 1$ we
    have   $\displaystyle   \lim_{\nu   \to   \frac{d(\beta-1)}{\beta}}
    N_\beta(\widetilde{A}) = 0$ while  $N_\alpha(A)$ remains finite for any $\nu
    >  0$:  $N_\alpha(A)   N_\beta(\widetilde{A})$  can  be  arbitrarily  small.
    Likewise,    when    $\beta<1$,     $\displaystyle    \lim_{\nu    \to    0}
    N_\beta(\widetilde{A}) = 0$ and  $\displaystyle \lim_{\nu \to 0} N_\alpha(A)
    =   +   \infty$,   but   $\displaystyle   \lim_{\nu   \to   0}   N_\alpha(A)
    N_\beta(\widetilde{A}) = 0$.}
\label{NaNb_nu_D0:fig}
\end{figure}

\begin{figure}[htbp]
\psfrag{nu}{$\nu$}
\psfrag{NaNbP}{\hspace{-9mm}\footnotesize $N_\alpha, N_\beta$ and $N_\alpha N_\beta$}
  \centerline{\includegraphics[width=6.5cm]{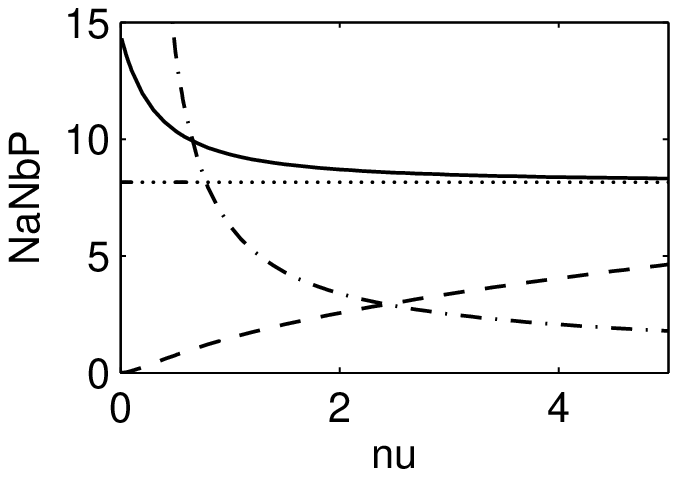}
    \hspace{5mm} \includegraphics[width=6.5cm]{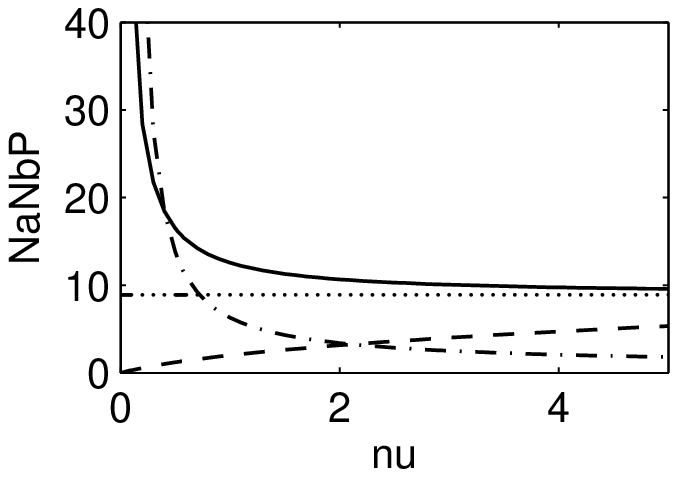}}
  \caption{Behavior  of $N_\alpha(A)$  (dash-dotted  line), $N_\beta(\tilde{A})$
    (dashed line),  and of  the product $N_\alpha(A)  N_\beta(\tilde{A})$ (solid
    line)  versus   $\nu$  in  the   Student-$t$  case,  to   illustrate  result
    (\ref{generalUP_C:eq}) (see  text).  Here $(\alpha,\beta) =  (2,2/3) \in \C$
    (left) and $(\alpha,\beta) = (2,1/2)  \in \D \backslash \C$ (right). In both
    cases one has , $\displaystyle  \lim_{\nu \to 0} N_\beta(\widetilde{A}) = 0$
    and $\displaystyle \lim_{\nu \to 0} N_\alpha(A) = + \infty$, but on $\C$ the
    product $  N_\alpha(A) N_\beta(\widetilde{A})$ has  a finite limit  while on
    $\D \backslash \C$  the limit is infinite. The  dotted line represents bound
    (\ref{generalBound_C:eq})         corresponding         to        inequality
    (\ref{generalUP_C:eq}).}
\label{NaNb_nu_CD:fig}
\end{figure}

From the  domain of existence  of $N_\beta(\widetilde{A})$, we  will distinguish
the cases $\beta >  1$, $\frac{1}{2} \le \beta < 1$ and  the limit case $\beta =
1$.  Playing  with the degree of  freedom $\nu$, we will  show that $N_\alpha(A)
N_\beta(\widetilde{A})$ can be arbitrarily small.
\begin{itemize}
\item Consider  first the case $\beta  > 1$. Then fix  $(\alpha,\beta) \in \D_0$
  and consider further  the case where $\frac{d (\beta-1)}{\beta}  < \nu \le d$.
  One  can easily check  that $\frac{\alpha  (d+\nu)-d}{2} >  \frac{d}{2} \left(
    \frac{\alpha}{\widetilde{\beta}}-1\right)$:        from       (\ref{NSt:eq})
  $N_\alpha(A)$   exists  and   is  finite   whatever  $\nu   \in  \left[\frac{d
      (\beta-1)}{\beta} \: , \: d \right]$.  Moreover, from the integral term of
  (\ref{NSttilde:eq})   and   \cite[9.6.9]{AbrSte70},   one   can   check   that
  $\displaystyle \lim_{\nu \to \frac{d (\beta-1)}{\beta}} N_\beta(\widetilde{A})
  =   0$.   As   a  consequence,   for   any  $(\alpha,\beta)   \in  \D_0   \cap
  \{(\alpha,\beta) | \beta > 1\}$,  we have $\displaystyle \lim_{\nu \to \frac{d
      (\beta-1)}{\beta}}  N_\alpha(A) N_\beta(\widetilde{A})  = 0$,  which proves
  that the product $N_\alpha(A) N_\beta(\widetilde{A})$ can be arbitrarily small
  in $\D_0 \cap \{ (\alpha,\beta) | \beta > 1\}$.
\item consider now  the case $ \frac{1}{2} \le  \beta < 1$, and $0  < \nu$. Then
  fix  again  a  pair  $(\alpha  ,   \beta)  \in  \D_0$.   One  can  check  that
  $\frac{\alpha  (d+\nu)-d}{2}   >  \frac{d}{2}  (\alpha-1)  >   0$.   The  last
  inequality comes  from $(\alpha ,  \beta) \in \D_0$,  in which if $\beta  < 1$
  then  $\alpha >  1$.  Hence,  from (\ref{NSt:eq})-(\ref{NSttilde:eq})  one can
  write
\begin{eqnarray*}
N_\alpha(A) N_\beta(\widetilde{A}) \propto \Gamma^{\frac{\alpha}{d(\alpha-1)} +
\frac{\beta}{d(\beta-1)}} \left( \frac{\nu}{2} \right)
\end{eqnarray*}
where the coefficient  of proportionality exists and is finite  for any $\nu \ge
0$. Since  in $\D_0$ one has $\frac{\alpha}{\alpha-1}  + \frac{\beta}{\beta-1} <
0$,   we    obtain   that   $\displaystyle   \lim_{\nu    \to   0}   N_\alpha(A)
N_\beta(\widetilde{A})  =   0$,  which  proves  that   the  product  $N_\alpha(A)
N_\beta(\widetilde{A})$ can be arbitrarily  small in $\D_0 \cap \{(\alpha,\beta)
| \beta < 1\}$.
\item The case $\beta=1$ can be deduced by continuity of $N_\beta$ as a function
  of $\beta$.
\end{itemize}
As a  conclusion, this example  shows that for  any pair of entropic  indices in
$\D_0$,  the  positive  quantity  $N_\alpha(A)  N_\beta(\widetilde{A})$  can  be
arbitrarily small.   This is  illustrated on figures  \ref{NaNb_nu_D0:fig} where
the  behavior  of  $N_\alpha(A)  N_\beta(\widetilde{A})$  versus  $\nu$  in  the
Student-$t$ case is  depicted, for chosen pairs $(\alpha,\beta)  \in \D_0$.  For
the  sake  of  comparison, let  us  recall  that  a  non-trivial bound  for  the
power-entropies product exists in the  case of pairs of entropic indices located
in  $\D$, as  given  by  Result \ref{generalUP_C:eq}.   This  is illustrated  in
Fig.~\ref{NaNb_nu_CD:fig},     where     the     behavior    of     $N_\alpha(A)
N_\beta(\widetilde{A})$ versus  $\nu$ in the  Student-$t$ case is  depicted, for
chosen pairs $(\alpha,\beta) \in \D$.

\ack
This work was performed thanks to the support of a CNRS/CONICET cooperation grant.

The authors  would like to acknowledge  the two anonymous  reviewers for helpful
comments on the original manuscript.


\bibliography{arXiv0709_3011_ZPV08}

\begin{thebibliography}{10}
\expandafter\ifx\csname url\endcsname\relax
  \def\url#1{\texttt{#1}}\fi
\expandafter\ifx\csname urlprefix\endcsname\relax\def\urlprefix{URL }\fi

\bibitem{Bia06}
I.~Bialynicki-{B}irula, Formulation of the uncertainty relations in terms of
  the {R\'e}nyi entropies, Physical Review A 74~(5) (2006) 052101.

\bibitem{ZozVig07}
S.~Zozor, C.~Vignat, On classes of non-{G}ausian asymptotic minimizers in
  entropic uncertainty principles, Physica A 375~(2) (2007) 499--517.

\bibitem{Lui07}
A.~Luis, Quantum properties of exponential states, Physical Review A 75 (2007)
  052115.

\bibitem{Deu83}
D.~Deutsch, Uncertainty in quantum measurements, Physical Review Letters 50~(9)
  (1983) 631--633.

\bibitem{PorPla96}
M.~Portesi, A.~Plastino, Generalized entropy as measure of quantum uncertainty,
  Physica A 225~(3-4) (1996) 412--430.

\bibitem{BiaMad85}
I.~Bialynicki-{B}irula, J.~L. Madajczyk, Entropic uncertainty relations for
  angular distributions, Physics Letters 108A~(8) (1985) 384--386.

\bibitem{CovTho91}
T.~M. Cover, J.~A. Thomas, Elements of Information Theory, John Wiley \& Sons,
  New-York, 1991.

\bibitem{MaaUff88}
H.~Maassen, J.~B.~M. Uffink, Generalized entropic uncertainty relations,
  Physical Review Letters 60~(12) (1988) 1103--1106.

\bibitem{HarLit52}
G.~Hardy, J.~E. Littlewood, G.~P{\'o}lya, Inequalities, 2nd Edition, Cambridge
  University Press, Cambridge, UK, 1952.

\bibitem{CosHer03}
J.~A. Costa, A.~O. {Hero III\@}, C.~Vignat, On solutions to multivariate
  maximum {$\alpha$}-entropy problems, in: A.~Rangarajan, M.~A.~T. Figueiredo,
  J.~Zerubia (Eds.), 4th International Workshop on Energy Minimization Methods
  in Computer Vision and Pattern Recognition (EMMCVPR), Vol. 2683 of Lecture
  Notes in Computer Sciences, Springer Verlag, Lisbon, Portugal, 2003, pp.
  211--226.

\bibitem{VigHer04}
C.~Vignat, A.~O. {Hero III\@}, J.~A. Costa, About closedness by convolution of
  the {T}sallis maximizers, Physica A 340~(1-3) (2004) 147--152.

\bibitem{ReeSim75}
M.~Reed, B.~Simon, Methods of Modern Mathematical Physics, vol II: Fourier
  Analysis, Self-Adjointness, Academic Press Inc., 1975.

\bibitem{Bec75:02}
W.~Beckner, Inequalities in {F}ourier analysis on ${R}^n$, Proceeding of the
  National Academy of Sciences of the USA 72~(2) (1975) 638--641.

\bibitem{BiaMyc75}
I.~Bialynicki-{B}irula, J.~Mycielski, Uncertainty relations for information
  entropy in wave mechanics, Communications in Mathematical Physics 44~(2)
  (1975) 129--132.

\bibitem{GadBen85}
S.~R. Gadre, R.~D. Bendale, Maximization of atomic information-entropy sum in
  configuration and momentum spaces, International journal of quantum chemistry
  28~(2) (1985) 311--314.

\bibitem{Lar90}
U.~Larsen, Superspace geometry: the exact uncertainty relationship between
  complementary aspects, Journal of Physics A 23~(7) (1990) 1041--1061.

\bibitem{ChaMon05}
K.~C. Chatzisavvas, C.~C. Moustakidis, C.~P. Panos, Information entropy,
  information distances, and complexity in atoms, The Journal of Chemical
  Physics 123~(17) (2005) 174111.

\bibitem{SenPan07}
K.~D. Sen, C.~P. Panos, K.~C. Chatzisavvas, C.~C. Moustakidis, Net {F}isher
  information measure versus ionization potential and dipole polarizability in
  atoms, Physics Letters A 364~(3-4) (2007) 286--290.

\bibitem{Lor54}
R.~Lord, The use of the {H}ankel transform in statistics {I}. {G}eneral theory
  and examples, Biometrika 41~(1/2) (1954) 44--55.

\bibitem{GraRyz80}
I.~S. Gradshteyn, I.~M. Ryzhik, Table of Integrals, Series, and Products,
  Academic Press, San Diego, 1980.

\bibitem{AbrSte70}
M.~Abramowitz, I.~A. Stegun, Handbook of Mathematical Functions with Formulas,
  Graphs, and Mathematical Tables, 9th printing, Dover, New-York, 1970.

\end{thebibliography}
\bibliographystyle{elsart-num}

\end{document}